# Distributed Finite-Time Termination of Consensus in the Presence of Delays

Mangal Prakash*, Saurav Talukdar*, Sandeep Attree, Vikas Yadav, *a*nd Murti V. Salapaka
*Both authors have contributed equally

*Abstract*—Linear consensus iterations guarantee asymptotic convergence, thereby, limiting their applicability in applications where consensus value needs to be used in real time to perform a system level task. It also leads to wastage of power and communication resources. In this article, an algorithm is proposed which enables each node to detect in a distributed manner and in finite number of iterations, when every agent in the network is within a user specified threshold of the consensus value (approximate consensus) and hence terminate further communications and computations associated with consensus iterations. This article develops a distributed algorithm for achieving this approximate consensus in presence of random time-varying bounded communication delays. Moreover, the article instantiates the algorithm developed to distributively determine the average of the initial values held by agents in finite number of iterations. Specifically, this algorithm relies on distributively determining the maximum and minimum of values held by the agents. The approach presented here offers several advantages, including reduced computational complexity, and hence, is suited for hardware implementation. An experimental test bed of Raspberry-Pi agents that communicate wirelessly over neighborhoods is employed as a platform to demonstrate the effectiveness of the developed algorithm.

## 1. Introduction

In recent times networks have gained widespread adoption for representing and analyzing large scale systems. Applications of the networks framework span multiple disciplines that include economics, neuroscience and social sciences [1] and emerging science and technology, including, "internet of things" [2]. Multi agent systems whose dynamics are governed by a network topology, often collaborate with each other in order to achieve system level objectives. Increasingly in applications, the size of the system imposes severe restrictions on resources that include computational capacity as well as communication bandwidth [3]. Due to these limitations coordination of multiple agents in large scale networked system calls for distributed algorithms.

A problem that has received considerable attention in co-ordination of multi-agent systems addresses the consensus problem, where how agents can compute a common value determined by initial values held by agents in a distributed manner is devised and analyzed [3], [4]. Here, all agents in the network strive to attain a common value of interest by sharing information with their neighbors in the network, which defines the communication layer. Consensus has found applications in parallel computers [5], distributed coordination of mobile autonomous agents [6], distributed data fusion in sensor networks [7], large scale power networks [8] and many more. A special case of consensus is that of average consensus where agents converge to the average of the initial conditions held by agents through local communication [4], [9]. There is considerable research toward algorithms for distributively reaching consensus with contributions from communications, control and computer science areas [10].

The convergence to the consensus value, in presence or absence of delays, when linear strategies are used is typically achieved asymptotically [11], [12], [13], [14], [15]. Here, agents keep updating their states and communicating with their neighbors forever leading to wastage of power and computational resources, which is untenable in resource constrained applications such as sensor networks where the resources available for each agent is limited. Moreover, in many real time applications including the power grid, distributed finite-time termination of consensus iterations is essential as the consensus value when determined is used in real-time by local systems to perform important tasks [16]. In such situations, it is necessary for each node to detect in finite-time, if approximate consensus is achieved within a specified error margin and thus, terminate computation as well as communication.

Distributed termination of the average consensus in finite time in the presence of time-varying but bounded delays is presented in [17]. [17] relies on computing and storing Hankel matrices at every iteration, which can be computationally expensive if the size of matrix is large; here the size depends on the network size as well as the number of iterations needed for convergence. We introduce a computationally efficient approach for finite time termination of consensus and average consensus both. We use maximum and minimum consensus to distributively determine the proximity (within a prespecified tolerance) of each node to the consensus/ average consensus value. Numerous applications of maximum and minimum consensus in other areas are reported in literature. [18] uses maximum consensus for synchronization of wireless sensor networks. Zhang and Li demonstrated the application of minimum consensus for tackling shortest path planning problem in graphs [19]. However, we explore a niche application of using the maximum and minimum consensus in an innovative manner to arrive at a distributed finite time termination criterion of consensus/ average consensus in the presence of fixed as well as time-varying link delays.

M. Prakash, S. Attree and M. V. Salapaka are with the Department of Electrical and Computer Engineering, University of Minnesota, Minneapolis, MN, 55455, USA e-mail: {praka027,murtis}@umn.edu.
S. Talukdar is with the Department of Mechanical Engineering, University of Minnesota, Minneapolis, MN, 55455, USA e-mail: taluk005@umn.edu.

*Contributions:* This article develops and analyzes an algorithm, where agents exchange information with other agents in their neighborhoods in the presence of unknown, time-varying but uniformly bounded communications delays, to achieve approximate consensus/average consensus (i.e, a termination criterion which can be used by agents to terminate computations distributively while realizing the consensus/ average consensus value is within the prespecified error tolerance). The algorithm depends on iteratively computing the maximum and minimum of the values held by agents with a small computational footprint. It is worth noting that the presentation in the manuscript is focused on the case of fixed communication delays to keep the discussion simple and the extension for the case of time varying communication delays is presented in the Appendix. Furthermore, the algorithm developed is instantiated to the average consensus problem as well and explicit bounds on the error from the average value due to finite-time termination of the average consensus protocol are determined. To the best of the knowledge of authors, the distributed finite-time termination algorithms developed here are simpler and computationally less expensive than other existing algorithms till date. The performance of the algorithms is illustrated using prototype networks of Raspberry Pis, where agents experience uncertainty in the communication channels. The algorithms proposed in this article for distributed finite time termination of consensus/ average consensus are generalizations and non-trivial extensions of the delay free framework presented in [20], which is the conference proceedings article by the authors. Notably, the algorithm presented in [20] fails when delays are present on communication channels and hence, the necessity arises for implementing the algorithm presented in this article.

The rest of the paper is organized as follows. In Section 2, the dynamics of the consensus algorithm executed by each agent is presented along with the needed notations, definitions and assumptions. In Section 3, the distributed finite-time stopping criterion for consensus algorithm based on maximum and minimum consensus protocols is developed. The average consensus protocol in presence of delays is presented in Section 4, followed by the development of distributed finite-time stopping criterion for average consensus protocol. The performance of the proposed distributed finite-time algorithms illustrated through simulations as well as with experiments on real communication networks realized through Raspberry Pi devices is presented in Section 5. Finally, the conclusions are presented in Section 6.

## 2. THE CONSENSUS PROTOCOL: BACKGROUND, DEFINITION AND ASSUMPTIONS

In this section definitions and notations needed for subsequent development (for details refer [21] and [22]) are provided. Consider the following definitions:

- *Directed and Undirected Graph*: A directed graph $G$ is a pair $\{V, E\}$ where $V$ is a set of vertices or nodes and $E$ is a set of edges, which are ordered subsets of two distinct elements of $V$. If an edge from $j \in V$ to $i \in V$ exists then it is denoted as $(i, j) \in E$. An undirected graph $G$ is a pair $\{V, E\}$ where $V$ is a set of vertices or nodes and $E$ is a set of edges such that for every pair of distinct nodes $i \in V$ and $j \in V$, if $(i, j) \in E$ then $(j, i) \in E$.
- *Directed Path*: In a directed graph, a directed path from node $j$ to $i$ exists if there is a sequence of distinct directed edges of $G$ of the form $(i, k_1), (k_1, k_2), ..., (k_m, j)$.
- *Strongly Connected Graph*: A directed graph is strongly connected if it has a directed path between each pair of distinct nodes $i$ and $j$.
- *In-neighbor of node*: Given a graph $G = \{V, E\}$ (directed or undirected), a node $j \in V$ is said to be an in-neighbor of node $i \in V$ if $(i, j) \in E$. The set of in-neighbors of node $i \in V$ is denoted by $N_i^- := \{j : (i, j) \in E\}$.
- *Diameter of Graph*: The longest shortest distance between any pair of nodes in a graph is the diameter of the graph and is denoted by $D$.
- *Stochastic Matrices*: A real $n \times n$ matrix $A = [a_{ij}]$ is called a row stochastic matrix if $1 \geq a_{ij} \geq 0$ for $1 \leq i, j \leq n$ and $\sum_{j=1}^{n} a_{ij} = 1$ for $1 \leq i \leq n$. A real $n \times n$ matrix $A = [a_{ij}]$ is called a column stochastic matrix if $1 \geq a_{ij} \geq 0$ for $1 \leq i, j \leq n$ and $\sum_{i=1}^{n} a_{ij} = 1$ for $1 \leq j \leq n$. A real $n \times n$ matrix $A = [a_{ij}]$ is called a doubly stochastic matrix if it is both row stochastic and column stochastic.
- *Non-negative Matrix and Primitive Matrix*: A non-negative matrix is a matrix all of whose entries are greater than or equal to zero. A primitive matrix is a square nonnegative matrix if it is irreducible and has only one eigenvalue of maximum modulus which is positive.

Consider a directed graph $G = \{V, E\}$. For every $(i, j) \in E$, a weight $p_{ij} > 0$ is associated which represents the weight, node $i$ gives to any information received from node $j$. $P = [p_{ij}]$ represents the weight matrix associated with graph $G$. If $(i, j) \notin E$ then $p_{ij} = 0$. In this article, the problem of consensus and average consensus on a network of agents in presence of delays on communication links is studied. The following section assumes that delays on the links are fixed and the extension of the results for the time varying but bounded communication delay case can be found in Appendix B. The consensus problem under consideration here admits the following assumptions.

A1. Graph $G$ is strongly connected.
A2. Weight matrix $P$ associated with the graph $G$ is row stochastic.
A3. Each edge in $G$ has a fixed delay, that is, for any two nodes $i, j \in V$ such that $(i, j) \in E$, the delay on the link $(i, j)$ is $\tau_{ij}(k) = \tau_{ij}$ for all time instants $k$.
A4. $(i, i) \in E$ and $\tau_{ii} = 0$ for all $i \in V$.
A5. The delay in the network is bounded and finite, that is, $\tau_{ij} \leq \bar{\tau}$, $\bar{\tau} \in \mathbb{R}$, for all $i, j \in V$.

The nodes do not have the knowledge of the delay associated with the edges in their in-neighborhood. Now, the consensus update rule for a directed network of agents interacting according to a fixed communication topology in the presence of fixed communication delays on the communication channels is presented. Under the assumptions A1-A5, agent $i$ updates its state $x_i(k+1)$ at $(k+1)^{th}$ iteration by taking a weighted linear combination of its own value and possibly delayed values of its in-neighbors received at $k^{th}$ iteration. The update rule is described by:

$$x_i(k+1) = p_{ii} x_i(k) + \sum_{j \in N_i^-} p_{ij} x_j(k - \tau_{ij}), \quad (1)$$

where, $N_i^-$ is determined by the graph $G = \{V, E\}$. Let $x(k)$ denote the column vector of all nodal states. Based on the update rule (1), consensus is defined as follows.

**Definition 1.** *(Consensus): A system of $n$ agents is said to have achieved consensus if for any set of initial conditions $x(0) \in \mathbb{R}^n$, there exists $\alpha \in \mathbb{R}$ such that $\lim_{k \to \infty} x_i(k) = \alpha$ for $i = 1, 2, ..., n$.*

Note that $\alpha$ does not depend on the node index $i$. It should be noted from (1) that at any update iteration $k$, from each of its in-neighbors node $i$ will receive one packet of information; albeit lagged information from a neighbor if there is a delay in the communication channel.

**Theorem 2.1.** [12], [13] *Update rule given by (1) under the assumptions mentioned in Section 2 achieves consensus.*

## 3. MAX-MIN CONSENSUS ALGORITHMS

In this section, first results based on the consensus update algorithm given by (1) are established and then Maximum and Minimum consensus algorithms are defined and their convergence is established. Subsequently a discussion on determining a distributed finite-time stopping criterion to detect if consensus is reached is provided.

Let for node $i$, the maximum over all values held by its neighbors including itself currently and in the $\bar{\tau}$ past instants be given by,

$$\{q_i, \tau_{q_i}\} := \arg\max_{\substack{j \in N_i^- \cup \{i\} \\ r=\{0,1,2,...,\bar{\tau}\}}} x_j(k - r), \quad (2)$$

for some $\tau_{q_i} \in \{0, 1, 2, ..., \bar{\tau}\}, q_i \in N_i^- \cup \{i\}$. Thus, $x_{q_i}(k - \tau_{q_i})$ is the maximum value in the neighborhood of node $i$ in the horizon of $\bar{\tau}$ into the past starting from iteration $k$. Similarly, let for node $i$, the minimum over all values held by its neighbors including itself, currently and in the $\bar{\tau}$ past instants be given by,

$$\{s_i, \tau_{s_i}\} := \arg\min_{\substack{j \in N_i^- \cup \{i\} \\ r=\{0,1,2,...,\bar{\tau}\}}} x_j(k - r), \quad (3)$$

for some $\tau_{s_i} \in \{0, 1, 2, ..., \bar{\tau}\}, s_i \in N_i^- \cup \{i\}$. Thus, $x_{s_i}(k - \tau_{s_i})$ is the minimum value in the neighborhood of node $i$ in the horizon of $\bar{\tau}$ into the past starting from iteration $k$. Furthermore, consider the maximum and minimum over all nodal values over the horizon $\{k - \bar{\tau}, ..., k-1, k\}$ in the past as given by,

$$M(k) := \max_{j \in V} x_{q_j}(k - \tau_{q_j}), \quad (4)$$

and,

$$m(k) := \min_{s \in V} x_{q_s}(k - \tau_{q_s}). \quad (5)$$

**Lemma 3.1.** *Consider the update rule (1). Then for all time instants $k' \geq k$ and for all $i \in V$,*

$$x_i(k') \leq \max_{j \in V} x_{q_j}(k - \tau_{q_j}) = M(k), \text{ and,} \quad (6)$$

$$x_i(k') \geq \min_{s \in V} x_{q_s}(k - \tau_{q_s}) = m(k). \quad (7)$$

*Proof.* See Appendix A for proof. □

*Lemma 3.1* establishes that the value held by an agent in the future is always bounded above by the maximum over the current and delayed values over a horizon $\bar{\tau}$ of all the nodal states. Moreover, the value held by the agent is bounded below by the minimum over the current and delayed values over a horizon $\bar{\tau}$ of all nodal states.

**Lemma 3.2.** *Consider a strongly connected graph $G = \{V, E\}$ running consensus protocol given by (1) with an initial condition $x(k)$. Let $i$ be a node such that $x_i(k) < M(k)$ and let $j$ be a node such that $x_j(k) > m(k)$, then for all time instants $k' \geq k$, $x_i(k') < M(k)$ and $x_j(k') > m(k)$.*

*Proof.* See Appendix A for proof. □

*Lemma 3.2* establishes that if the value held by an agent $i$ at the present instant of time is strictly less (greater) than the maximum (minimum) over the current and delayed values over a horizon $\bar{\tau}$ of all the nodal states, then the value of agent $i$ continues to be strictly less (greater) than this maximum (minimum) for all future instants.

Consider the maximum and minimum value in the network, which is defined as, $\max x(k) := \max_{i \in V} x_i(k)$ and $\min x(k) := \min_{i \in V} x_i(k)$ respectively.

**Lemma 3.3.** *Consider a strongly connected graph $G = \{V, E\}$ with an update rule for the consensus protocol given by (1) with an initial condition $x(k)$ such that $\min x(k) < \max x(k)$. Then for all $k' \geq k + D(\bar{\tau}+1), \max x(k') < M(k)$ and $\min x(k') > m(k)$.*

*Proof.* See Appendix A for proof. □

*Lemma 3.3* establishes that if the maximum value over all nodal states at the present instant is strictly greater than the minimum value over all nodal states at the present instant, then the maximum (minimum) over all nodal states after $D(1+\bar{\tau})$ instants in the future will be strictly less (greater) than the present maximum (minimum) over the current and delayed values over a horizon $\bar{\tau}$ of all the nodal states.

Define by $T^* := D(1+\bar{\tau}) + \bar{\tau}$. Note that $T^*$ is a constant for a fixed interconnection topology.

**Theorem 3.1.** *Consider a strongly connected graph $G = \{V, E\}$ with an update rule for the consensus protocol given by (1) and an initial condition $x(lT^*)$ such that $\min x(lT^*) < \max x(lT^*)$, where $l \geq 0$. Then, $M((l+1)T^*) < M(lT^*)$ and $m((l+1)T^*) > m(lT^*)$.*

*Proof.* Using $k = lT^*$ in Lemma 3.3, it follows that for $k' \geq lT^* + D(1 + \bar{\tau})$,

$$\max x(k') < M(lT^*). \tag{8}$$

By definition,

$$M((l+1)T^*) = \max_{j \in V} x_{q_j}((l+1)T^* - \tau_{q_j}).$$

Since the index $((l+1)T^* - \tau_{q_j}) \geq lT^* + D(1+\bar{\tau})$, it follows that

$$M((l+1)T^*) < M(lT^*).$$

The rest of the proof is left to the reader. □

Theorem 3.1 implies that $M(lT^*)$ is a strictly decreasing sequence and $m(lT^*)$ is a strictly increasing sequence as a function of the index $l$.

**Corollary 3.1.** *Consider a strongly connected graph $G = \{V, E\}$ running consensus protocol given by (1) with an initial condition $x(lT^*)$ such that $\min x(lT^*) < \max x(lT^*)$, where, $l \geq 0$. Then, (a) $\max x((l+1)T^*) < M(lT^*)$ and $\min x((l+1)T^*) > m(lT^*)$. Also, (b) $\max x((l+1)T^*) - \min x((l+1)T^*) < M(lT^*) - m(lT^*)$.*

*Proof.* The proof of *(a)* follows directly from Theorem 3.1. The proof of *(b)* is a direct consequence of *(a)*. □

**Corollary 3.2.** *Consider the consensus protocol given by (1) with each node converging to $\alpha$, i.e., $\lim_{k \to \infty} x_i(k) = \alpha$ for all $i = 1, 2, ...n$. Then the sequences $\{M(lT^*)\}_{l \in \mathbb{N}}$ and $\{m(lT^*)\}_{l \in \mathbb{N}}$ converge to $\alpha$ as $l \to \infty$. Further, the sequence $\{M(lT^*) - m(lT^*)\}_{l \in \mathbb{N}} \to 0$ as $l \to \infty$.*

*Proof.* From the hypothesis it follows that, there exist $\alpha$ such that,

$$\lim_{k \to \infty} x_i(k) = \alpha, \text{for all } i = 1, 2, ...n.$$

Further, $M(lT^*)$ and $m(lT^*)$ are subsequences of convergent sequence $x(k)$ and hence converge to the same limit $\alpha$. Thus both $M(lT^*)$ and $m(lT^*)$ converge to $\alpha$ as $l \to \infty$. This implies,

$$M(lT^*) - m(lT^*) \to 0 \text{ as } l \to \infty.$$

□

*Corollary 3.1 (b)* and *Corollary 3.2* together imply that $\max x((l+1)T^*) - \min x((l+1)T^*) \to 0$ as $l \to \infty$. In what follows an algorithm is devised which converges to $\max x((l+1)T^*)$ and $\min x((l+1)T^*)$ in finite number of iterations with the finite-time stopping criteria based on the difference between $\max x((l+1)T^*)$ and $\min x((l+1)T^*)$ and evaluating whether the difference is less than the user specified error tolerance. Towards this end, the maximum/ minimum consensus protocol are developed in the following subsection and the finite-time convergence of maximum/ minimum consensus protocols in presence of fixed communication delays is proved. Preceding results are utilized to develop a finite-time termination criterion for the consensus protocol.

Define $\tilde{T} := (\bar{\tau} + 1)$. Note that in the subseDefnt discussions, these two notations will be used interchangeably. ine the queindicator function for the link from node $j \in V$ to node $m \in V$ at time $k$ as

$$I_{k,mj}(\tau) = \begin{cases} 1, & \text{if } \tau_{mj}(k) = \tau \\ 0, & \text{if } \tau_{mj}(k) \neq \tau, \end{cases}$$

### A. Maximum and Minimum Consensus Protocols

The Maximum Consensus Protocol denoted by MXP computes the maximum of the given initial node conditions $z(0) = [z_1(0)\ z_2(0)....z_n(0)]^T$ in a distributed manner. It takes $z(0)$ as an input and generates a sequence of node values based on the following update rule for node $m$ for $k \geq 0$,

$$z_m(k\bar{\tau} + l) = z_m(k\bar{\tau} + l - 1), l \in \{k+1, \cdots, k+\bar{\tau}\}, \tag{9}$$

$$z_m((k+1)\tilde{T}) = \max_{j \in N_m^- \cup \{m\}} \{z_m((k+1)\tilde{T} - (r+1))I_{(k+1)\tilde{T}-r,mj}(r)\}_{r=0,1,...,\bar{\tau}}. \tag{10}$$

where, the indicator function is defined as in (9).

The Minimum Consensus Protocol denoted by MNP computes the minimum of the given initial node conditions $y(0) = [y_1(0)\ y_2(0)....y_n(0)]^T$ in a distributed manner. It takes $y(0)$ as an input and generates a sequence of node values $y(k)$ based on the following update rule for $k \geq 0$:

$$y_m(k\bar{\tau} + l) = y_m(k\bar{\tau} + l - 1), l \in \{k+1, \cdots, k+\bar{\tau}\}, \tag{11}$$

$$y_m((k+1)\tilde{T}) = \min_{j \in N_m^- \cup \{m\}} \{y_m((k+1)\tilde{T} - (r+1))I_{(k+1)\tilde{T}-r,mj}(r)\}_{r=0,1,...,\bar{\tau}}. \tag{12}$$

where, the indicator function is defined as in (9).

Note that (9) maintains value of $z_m$ at $z_m((k-1)\bar{\tau})$ till the $k^{th}$ epoch $k\bar{\tau}+l, l \in \{1, 2, ..., \bar{\tau}\}$ ends. On the other hand (10) updates $z_m$ at time instances which are multiples of $\tilde{T} = \bar{\tau} + 1$ based on recent information from the neighbors and itself. Effectively every $z_m$ update takes place once after every $\bar{\tau}$ iterations. Similarly, every $y_m$ update takes place once after every $\bar{\tau}$ iterations. MNP is similar to MXP since the minimum over a set of values is the negative of the maximum of the negative of the values. Next, results are established which will be useful to prove the convergence of MXP and MNP running through the update rules (9), (10), (11) and (12) in finite-time. Let $\tilde{z} := \max_{i \in V} z_i(0)$ and Let $\tilde{y} := \min_{i \in V} y_i(0)$.

**Lemma 3.4.** *Consider a directed graph $G = \{V, E\}$ with fixed delays with uniform bound $\bar{\tau}$ and an update rule for the Maximum Consensus Protocol (MXP) given by (9) and (10) and the Minimum Consensus Protocol (MNP) given by (11) and (12).*

(a) Then, for all $k' > 0$, $\max_{i \in V} z_i(k') = \tilde{z}$ and $\min_{i \in V} y_i(k') = \tilde{y}$.
(b) Let for some $k$, $z_j(k) = \tilde{z}$, that is, node $j$ has the maximum value in the network at the $k^{th}$ time instant. Then, for all instants $k' > k$, $z_j(k') = \tilde{z}$, that is node $j$ continues to be the maximum for $k' > k$.
(c) Let for some $k$, $y_{j'}(k) = \tilde{y}$, that is, node $j'$ has the minimum value in the network at the $k^{th}$ time instant. Then, for all instants $k' > k$, $y_{j'}(k) = \tilde{y}$, that is node $j'$ continues to be the minimum for $k' > k$.

*Proof.* (a) The proof is left to the reader. (b) The proof follows from the fact that if the maximum value at the current iteration is held at node $j$, then node $j$ continues to hold the maximum value in future iterations as well, as the update step of the MXP (10) includes the past value of node $j$. (c) The proof is similar to the proof of (b). □

**Lemma 3.5.** *Consider a directed graph $G = \{V, E\}$ with fixed delays with uniform bound $\bar{\tau}$ and an update rule for the Maximum Consensus Protocol (MXP) given by (9) and (10) and that for for the Minimum Consensus Protocol (MNP) given by (11) and (12). Let $z_{\pi_1}(0) = \max_{j \in V} z_j(0)$ and $y_{\pi'_1}(0) = \min_{j \in V} y_j(0)$. Then, for all $k \geq D(1 + \bar{\tau})$, and any $m \in V$,*

(a) $z_m(k) = z_{\pi_1}(0)$.
(b) $y_m(k) = y_{\pi'_1}(0)$.

*Proof.* (a) As $z_{\pi_1}(0) = \max_{j \in V} z_j(0)$ it follows from Lemma 3.4 that

$$z_{\pi_1}(k) = z_{\pi_1}(0) \text{ for all } k \geq 0.$$

Consider any node $i \in V$. Since the graph $G$ is strongly connected, there exists a directed path $(\pi_2, \pi_1)(\pi_3, \pi_2)...(i, \pi_d)$ connecting $i$ and $\pi_1$. It follows from the update rule (9) and (10) that within $\bar{\tau}$ iterations, $\pi_2$ will have received the value $z_{\pi_1}(0)$ and thus, $z_{\pi_2}(\bar{\tau}+1) = z_{\pi_1}(0)$; and for any $k \geq \bar{\tau}+1$, $z_{\pi_2}(k) = z_{\pi_1}(0)$. Using the above steps for $\pi_3, \pi_4, ..., \pi_d, i$, it follows that,

$$z_i(k) = z_{\pi_1}(0) \text{ for any } k \geq d(\bar{\tau}+1)$$

Thus if $k \geq D(\bar{\tau}+1) \geq d(\bar{\tau}+1)$; $z_i(k) = z_{\pi_1}(0)$. Since $D$ is the diameter of the graph $G$, it follows that, for $k \geq D(\bar{\tau}+1)$,

$$z_m(k) = z_{\pi_1}(0) = \max_{j \in V} z_j(0) \text{ for all } m \in V.$$

(b) The proof is similar to the proof of (a). This proves the theorem. □

### B. Distributed finite-time Algorithm for Terminating Consensus protocol

In this section, an algorithm based on Maximum-Minimum Consensus for stopping the consensus protocol distributively in finite-time is proposed using the results derived in the previous sub-sections. At time instants $T(j) = jT^*$, for $j = 1, 2, ...$ MXP and MNP protocols are reset with initial conditions $x(jT^*)$, thus $z(T(j)) = x(T(j))$ and $y(T(j)) = x(T(j))$.

Let,
$$\bar{\alpha}(j) := \max z(T(j)), \underline{\alpha}(j) := \min y(T(j)), \text{ and}$$
$$\beta(j) := \bar{\alpha}(j) - \underline{\alpha}(j).$$

**Lemma 3.6.** *Consider the consensus protocol given by (1) with each node converging to $\alpha$, i.e., $\lim_{k \to \infty} x_i(k) = \alpha$ for all $i = 1, 2, ...n$. Then the sequences $\bar{\alpha}(j)$ and $\underline{\alpha}(j)$ converge to $\alpha$ as $j \to \infty$. Further, the sequence $\beta(j) \to 0$ as $j \to \infty$.*

*Proof.* If $\{x(k)\}_{k=1, i \in V}^{\infty}$ converges to $\alpha$, that is, $\lim_{k \to \infty} x_i(k) = \alpha$ for all $i = 1, 2, ...n$, then it follows that $\lim_{k \to \infty} \max x_i(k) = \alpha$ and $\lim_{k \to \infty} \min x_i(k) = \alpha$ for all $i = 1, 2, ...n$.

Note that, $\bar{\alpha}(j)$ and $\underline{\alpha}(j)$ are subsequences of convergent sequences $\max x_i(k)$ and $\min x_i(k)$ respectively. Thus, $\bar{\alpha}(j)$ and $\underline{\alpha}(j)$ converge to the same limit $\alpha$. This implies that, $\beta(j) \to 0$ as $j \to \infty$. □

It should be noted that Lemma 3.6 is a consequence of Corollary 3.2.

---

**Algorithm 1:** Finite-time termination of consensus in presence of uniformly bounded delays

1 **Input:**
2    $x(0), D, \bar{\tau}, \rho, P$ ;         // Initial condition
3 **Initialize:**
4    $k := 0$;
5    $l := 1$;
6    $z_i := x_i(0)$;
7    $y_i := x_i(0)$;
8    $\theta := 1$;
9    $\psi := 0$;
10 **Repeat:**
11    $x_i(k+1) = p_{ii} x_i(k) + \Sigma_{j \in N_i^-} p_{ij} x_j(k - \tau_{ij})$
12    **if** $k + 1 = \psi + l(1 + \bar{\tau})$ **then**
       /* maximum and minimum consensus updates given by (10) and (12) for each node $i \in V$ */
13        $z_i := \max_{j \in N_i^- \cup \{i\}} z_i$;
14        $y_i := \min_{j \in N_i^- \cup \{i\}} y_i$;
15        $l := l + 1$
16    **end**
17    **emit:** $x_i(k+1)$, $y_i$ and $z_i$
18    **if** $k + 1 = \theta T^*$ **then**
19        **if** $z_i - y_i < \rho$ **then**
20            **break** ;     // stop $x_i$, $y_i$ and $z_i$ updates
21        **else**
22            $z_i := x_i(\theta T^*)$;
23            $y_i := x_i(\theta T^*)$;
24            $\theta := \theta + 1$;
25            $l := 1$ ;           // Reset
26            $\psi := \psi + k$
27        **end**
28    **end**
29    $k := k + 1$

**Theorem 3.2.** *Algorithm 1 converges in some finite-time $T_c < \infty$, that is, Algorithm 1 reaches line number 20. Mathematically, given any $\rho > 0$, there exists $T_c \in \mathbb{N}$ such that $\beta(j) < \rho$.*

*Proof.* It follows from Lemma 3.6 that $\beta(j) \to 0$ as $j \to \infty$. Thus, for any given $\rho > 0$, there exists an integer $j_0$ such that $\beta(j) < \rho$ for all $j \geq j_0$. This implies that the *Algorithm 1* converges in finite-time $T_c = T(j_0)$. □

**Remark 1.** *If $\bar{\tau} = 0$, Algorithm 1 and the results presented above reduce to those presented in [20].*

**Remark 2.** *All the results presented till this point for consensus and max-min consensus hold true even for uniformly bounded time varying delays under the assumption that only one packet of information from each neighbor can be processed at any time instant in the update rule (1). The asymptotic convergence of consensus under random delay model with at most one packet of information being processed at any time instant by any node is established in [12], [13] and [23]. The contribution here is that Algorithm 1 is valid even for these communication models for finite-time termination of consensus. The results and proofs for this communication model can be found in Appendix B.*

## 4. AVERAGE CONSENSUS PROTOCOL

Average consensus problem is a special case of consensus where all the nodes converge to the average of the initial conditions. and is defined as follows.

**Definition 2.** *(Average Consensus) A system of $n$ agents is said to have achieved average consensus if for any initial condition $x(0) \in \mathbb{R}^n$, $\lim_{k\to\infty} x_i(k) = \frac{\sum_{i=1}^n x_i(0)}{n}$ for all $i = 1, 2, ..., n$.*

In [24], it is shown that average consensus can be reached by using the ratio of two consensus updates as described below. First the assumptions and algorithm developed in [24] are discussed. Consider a directed graph with $n$ nodes which satisfies the following assumptions.

B1. Weight matrix $P$ associated with the directed graph is primitive and column stochastic.
B2. The directed graph is strongly connected.
B3. Any node $i$ in the directed graph has access to its own value at any instant $k$ without any delay.
B4. The delay on the directed edge connecting any two nodes $i$ and $j$ in the directed graph is bounded by some constant $\bar{\tau}$, i.e., $\tau_{ij} \leq \bar{\tau} < \infty$.

**Theorem 4.1.** *([24]) Suppose the assumptions B1-B4 are satisfied. Let $x_i(k)$ and $w_i(k)$ be the result of iterations*

$$x_i(k+1) = p_{ii}x_i(k) + \sum_{j \in N_i^-}\sum_{r=0}^{\bar{\tau}} p_{ij}x_j(k-r)I_{ij}(r), \text{ and} \quad (13)$$

$$w_i(k+1) = p_{ii}w_i(k) + \sum_{j \in N_i^-}\sum_{r=0}^{\bar{\tau}} p_{ij}w_j(k-r)I_{ij}(r). \quad (14)$$

*Let the initial conditions be given by $x(0) = [x_1(0)\ x_2(0)...x_n(0)]^T$ and $w(0) = 1_n$ where $1_n$ is a $n \times 1$ column vector of all ones. Then the ratio of $x_i(k)$ and $w_i(k)$ asymptotically converges to $\lim_{k\to\infty} \mu_j(k) = \frac{\sum_{i=1}^n x_i(0)}{n}$ for all $j = 1, ..., n$ where $\mu_j(k) := x_j(k)/w_j(k)$.*

In order to satisfy the column stochastic assumption of Theorem 4.1 and to extend Algorithm 1 (which requires row stochasticity) for average consensus, doubly stochastic weight matrices are chosen. A square matrix is irreducible if and only if its associated graph is strongly connected [25]. Using Perron Frobenius Theorem, it follows that doubly stochastic weight matrix $P$ is primitive [22]. Using Theorem 4.1 it can be shown that running two consensus protocols given by (13) and (14), the average consensus can be asymptotically achieved with the initial conditions as $x(0) = [x_1(0)\ x_2(0)...x_n(0)]^T$ and $w(0) = 1_n$ where $1_n$ is a $n \times 1$ column vector of all ones. Thus the ratio of $x_j(k)$ and $w_j(k)$ asymptotically converges to $\frac{\sum_{i=1}^n x_i(0)}{n} = c$, for all $j = 1, ..., n$.

### A. Distributed finite-time Algorithm for terminating Average Consensus Protocol

An MXP and an MNP associated with (13) are executed. Another MXP and MNP associated with (14) are also executed. By Theorem 3.1, both (13) and (14) converge. Let $\lim_{k\to\infty} x_i(k) = \alpha$ for all $i = 1, 2, ..., n$ and $\lim_{k\to\infty} w_i(k) = \sigma$ for all $i = 1, 2, ..., n$. Using Theorem 3.2, the consensus protocol given by (13) can be stopped in some finite-time $T_{c1}$ when $|x_i(T_{c1}) - \alpha| < \rho$. Also, using Theorem 3.2, the consensus protocol given by (14) can be stopped in some finite-time $T_{c2}$ when $|w_i(T_{c2}) - \sigma| < \rho$. Using Theorem 4.1 given $\epsilon > 0$, (the bound within which the deviation of states from average of initial conditions is permitted), $\rho$ can be chosen such that stopping the consensus protocols given by (13) and (14) in finite-time depending on chosen $\rho$, will ensure that average consensus can be achieved within a positive constant $\epsilon$ of the average of initial conditions $c$ in finite-time. The deviation from average of initial conditions by stopping the consensus protocols (13) and (14) in finite-time is now quantified in the following discussion.

It is shown in [24] that

$$\lim_{k\to\infty} \frac{x_i(k)}{w_i(k)} = c = \frac{\alpha}{\sigma}. \quad (15)$$

Since, the two consensus protocols given by (13) and (14) terminate when they reach within some specified bound $\rho > 0$, $x_i(k)$ and $w_i(k)$ may not converge to $\alpha$ and $\sigma$ respectively but instead attain $\tilde{\alpha}$ and $\tilde{\sigma}$ respectively as the terminal values. The ratio $\frac{x_i}{w_i}$ will deviate from $c = \frac{\alpha}{\sigma}$ depending on when the consensus protocols given by (13) and (14) are stopped. The deviation from average can be quantified explicitly as a function of $\rho, \tilde{\alpha}$ and $\tilde{\sigma}$, all of which are known values once $\rho$ is specified. The deviation is quantified by (16) below.

$$\frac{\tilde{\alpha}}{\tilde{\sigma}} - \frac{\tilde{\alpha}}{\tilde{\sigma}}\left(\frac{1 - \frac{\rho}{\tilde{\alpha}}}{1 + \frac{\rho}{\tilde{\sigma}}}\right) \geq \frac{x_i}{w_i} - c \geq \frac{\tilde{\alpha}}{\tilde{\sigma}} - \frac{\tilde{\alpha}}{\tilde{\sigma}}\left(\frac{1 + \frac{\rho}{\tilde{\alpha}}}{1 - \frac{\rho}{\tilde{\sigma}}}\right) \quad (16)$$

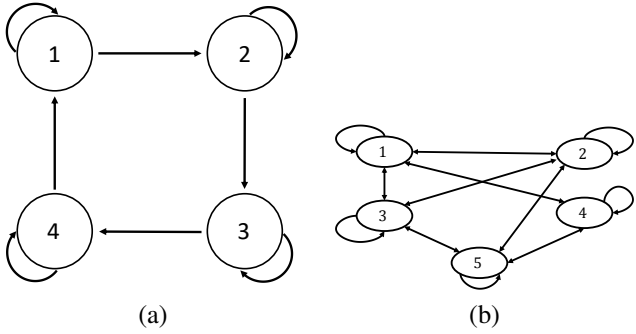

Figure 1: (a) Directed ring network of 4 nodes (b) A general representative network of 5 nodes.

Given the error tolerance $\epsilon$ within which the convergence to average of initial conditions is desired, $\rho$ should be chosen to satisfy the bounds in (17) and (18),

$$\frac{x_i}{w_i} - c \leq \frac{\tilde{\alpha}}{\tilde{\sigma}} - \frac{\tilde{\alpha}}{\tilde{\sigma}}\left(\frac{1 - \frac{\rho}{\tilde{\alpha}}}{1 + \frac{\rho}{\tilde{\sigma}}}\right) \leq \epsilon, \qquad (17)$$

$$\frac{x_i}{w_i} - c \geq \frac{\tilde{\alpha}}{\tilde{\sigma}} - \frac{\tilde{\alpha}}{\tilde{\sigma}}\left(\frac{1 + \frac{\rho}{\tilde{\alpha}}}{1 - \frac{\rho}{\tilde{\sigma}}}\right) \geq -\epsilon. \qquad (18)$$

This implies that given the error tolerance $\rho$ for $x_i(k)$ and $w_i(k)$, $\left|\frac{x_i(k)}{w_i(k)} - c\right| \leq \frac{\tilde{\alpha}}{\tilde{\sigma}} - \frac{\tilde{\alpha}}{\tilde{\sigma}}\left(\frac{1+\frac{\rho}{\tilde{\alpha}}}{1-\frac{\rho}{\tilde{\sigma}}}\right)$, $\rho$ can be chosen appropriately to ensure that the error from the average is as small as desired. A detailed derivation of (16) is given in Appendix C.

Using the Algorithm 1, distributed finite-time termination algorithm for average consensus is presented next.

---

**Algorithm 2:** Finite-time termination of average consensus in presence of uniformly bounded delays

1. Given $\epsilon$ (the permitted deviation from average of initial conditions), Choose appropriately $\rho$ for stopping the consensus updates given by (13) and (14)
2. Run *Algorithm 1* for both (13) and (14) maintained by each node $i \in V$
3. Stop *Algorithm 1* for both (13) and (14) together, i.e. for any node *Algorithm 1* terminates only when both (13) and (14) have met the termination criterion $\rho$. This terminates *Algorithm 2*
4. Use (16) to check the deviation of the so computed value from the actual average of initial conditions. If the deviation is within desired bound $\epsilon$, stop. Else repeat *Steps 1-3* with a smaller value of $\rho$.

---

## 5. RESULTS AND DISCUSSION

### A. Simulation Results

Here results of *Algorithm 2* for finite-time termination of average consensus on the ring network shown in Figure 1(a) and a representative 5 node network shown in Figure 1(b) in presence of fixed communication delays is presented. The communication weight matrix for the ring network is chosen to be,

$$P = \begin{bmatrix} 1/2 & 1/2 & 0 & 0 \\ 0 & 1/2 & 1/2 & 0 \\ 0 & 0 & 1/2 & 1/2 \\ 1/2 & 0 & 0 & 1/2 \end{bmatrix},$$

and, for the 5 node network the weight matrix is chosen to be

$$P = \begin{bmatrix} 2/5 & 1/5 & 1/5 & 1/5 & 0 \\ 1/5 & 2/5 & 1/5 & 0 & 1/5 \\ 1/5 & 1/5 & 2/5 & 0 & 1/5 \\ 1/5 & 0 & 0 & 2/5 & 2/5 \\ 0 & 1/5 & 1/5 & 2/5 & 1/5 \end{bmatrix}.$$

For the ring network in Figure 1(a), the delay on the edges are

$$\tau = \begin{bmatrix} 0 & 1 & 0 & 0 \\ 0 & 0 & 2 & 0 \\ 0 & 0 & 0 & 1 \\ 2 & 0 & 0 & 0 \end{bmatrix}.$$

The initial conditions are set as $x(0) = [50\ 70\ 150\ 30]^T$. The error tolerance for deviation from average of initial conditions is set to $\epsilon = 0.1$. The stopping bound is set to $\rho = 0.01$ for both the consensus algorithms given by (13) and (14). For the network in Figure 1(a), the consensus algorithm given by (13) converges to the value $42.85$ in 166 iterations and the consensus algorithm given by (14) converges to the value $0.57$ in 166 iterations. Accordingly the ratio of (13) and (14) converges to $74.99$ in 166 iterations as shown in Figure 2. The average of the initial conditions is 75 and thus, the deviation from the average achieved by implementing finite-time termination algorithm on ratio consensus is well within the bounds given by (16), (17) and (18).

For the network in Figure 1(b) the delay on the edges are

$$\tau = \begin{bmatrix} 0 & 2 & 1 & 0 & 0 \\ 2 & 0 & 0 & 0 & 3 \\ 1 & 0 & 0 & 0 & 0 \\ 0 & 0 & 0 & 0 & 0 \\ 0 & 3 & 0 & 0 & 0 \end{bmatrix}.$$

The initial conditions for this network are set as $x(0) = [1000\ 0\ 200\ 100\ 700]^T$ and the error tolerance for deviation from average of initial conditions is again set to be $\epsilon = 0.1$. The consensus algorithm given by (13) converges to the value $270.3$ in 40 iterations. For this network, the consensus algorithm given by (14) converges to the value $0.6757$ in 40 iterations. Accordingly the ratio of (13) and (14) converges to $399.9974$ in 40 iterations as shown in Figure 3. The average of the initial conditions is 400 and thus, the deviation from the average achieved by implementing finite-time termination algorithm on ratio consensus is well within the bounds given by (16), (17) and (18).

Table I presents the results for comparison of *Algorithm 2* with the distributed finite-time algorithm proposed in [17]. The computational complexity is listed in terms of the total number of nodes in the network and it can be seen that the proposed algorithm is more efficient in terms of both CPU as well as memory requirements.

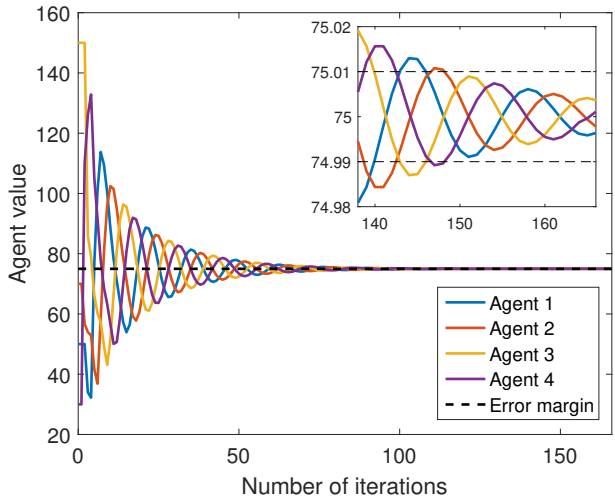

Figure 2: Simulation results for Maximum-Minimum consensus based distributed finite-time termination of average consensus for the ring network of 4 nodes shown in Figure 1(a).

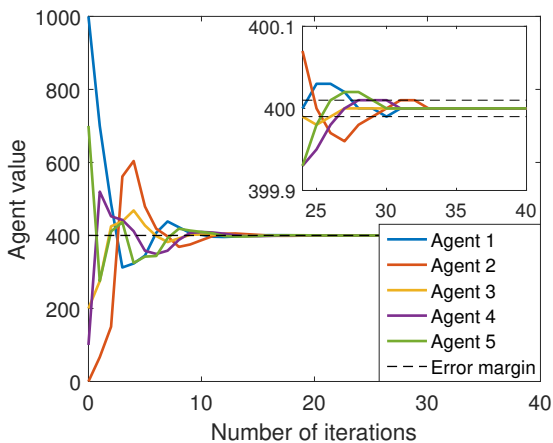

Figure 3: Simulation results for Maximum-Minimum consensus based distributed finite-time termination of average consensus for the network of 5 nodes shown in Figure 1(b).

Table I: Comparison of the simulation results for the 5 node network using the proposed algorithm and the algorithm in [17].

| Algorithm | Number of iterations | Computational complexity at a node for an iteration | |
|---|---|---|---|
| | | Time complexity | Space complexity |
| Algorithm in [17] | 30 | $\mathcal{O}(n^2)$ | $\mathcal{O}(n^2)$ |
| *Algorithm 2* | 48 | $\mathcal{O}(n)$ | $\mathcal{O}(n)$ |

Next, we demonstrate the applicability of our distributed finite time termination algorithm (*Algorithm 1*) for termination of consensus iterations in the presence of random (time varying) but bounded communication delays on the network shown in Figure 1(b). It is assumed that at each link at each time instant, the delay is a non-negative integer upper bounded by $\bar{\tau} = 3$. All the delays are sampled from a uniform distribution on $\{0,1,2,3\}$ for the simulation. The initial conditions are set as $x(0) = [1000\ 0\ 200\ 100\ 700]^T$. Each node determines that the consensus has reached within an error margin of $\rho = 0.01$ in 48 iterations. All the nodes converge to a common value of 261.1 (see Figure 4). It is worth noting that the final value at which the consensus protocol converges in presence of time-varying delays depends on the specific realization of the communication delays.

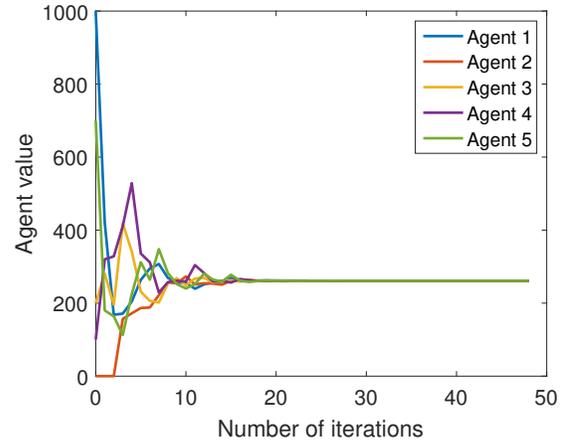

Figure 4: Simulation results for Maximum-Minimum consensus based distributed finite-time termination in the presence of uniformly bounded time-varying delays for the network of 5 nodes shown in Figure 1(b).

### B. Experimental Demonstration

In this section, the functionality and efficacy of the proposed algorithm on a physical network is established. First, the experimental setup is described and then the results and observations are discussed.

Rapsberry Pi devices [26] are used to setup the physical network for experimentation. Each of such devices is a Raspberry Pi 3 model b with configuration of 1.2 GHz CPU, 1 GB RAM and 802.11n Wireless LAN support [27]. Each of the agent nodes in the network is an individual Raspberry Pi device capable of communicating with other agents over Internet via a Wi-Fi network. In a practical setting, all nodes communicate with latency in the communication channel. The distribution of pair-wise communication latency for each node in the experimental realization of the network given in Figure 1(b) is shown in Figure 5. For instance, the box plot labelled as 'Agent 1' depicts the latency experienced by Raspberry Pi 1 for several communication requests sent to every other Raspberry Pi in the network i.e. Agents 2, 3, 4 and 5. It can be seen that the network latency ranges from 10ms to 100ms. Data communications happen via HTTP protocol [28], which is designed over TCP/IP [29] and ensures guaranteed, in-order delivery of data packets between the nodes. A NodeJS [30] based application capable of bi-directional communication for running the consensus algorithm on each Raspberry Pi is developed. NodeJS has an event-driven architecture capable of asynchronous I/O which greatly enhances throughput and scalability in real-time web applications making it an ideal choice for applications targeted by the framework of the article. The application accepts a configuration setup to initialize the consensus algorithm and communicate with other nodes for information exchange, a necessary requirement for a Raspberry Pi to act as an agent. All the agents are initialized by passing the configuration information before starting the consensus protocol, and the agent node will then transmit or receive data from the neighbouring nodes in real-time. The agent node also logs data periodically for analysis.

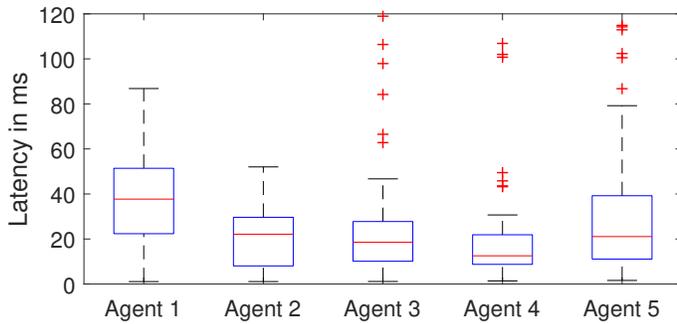

Figure 5: Distribution of pair-wise latency (in milliseconds) for each node in the network shown in Figure 1(b).

To illustrate the validity and performance of *Algorithm 2*, this algorithm is compared with the delay-free termination algorithm presented in [20] . Both the algorithms are tested in presence of constraints and variabilities that a physical network inherits. Figure 6 and Table II illustrate the performance of the delay-free termination algorithm for the 5 node network shown in Figure 1(b) and it can be seen that this algorithm does not terminate near the average of the nodal initial conditions. Next, the distributed finite-time termination of average consensus algorithm in presence of fixed delays is tested on the directed ring network having four agents as depicted in Figure 7(a) and the observed results are presented in Table III. The outcome of the same algorithm when applied to the network of 5 agents shown in Figure 1 (b) are presented in Figure 7(b) and Table IV. The delays and initial conditions for both the cases were set up in the same way as they were for simulations and the error margin i.e. $\rho$ is also chosen to be 0.01. From Tables III and IV and Figure 7 (a) and 7 (b), it is clear that the finite time termination algorithm terminates when the nodes reach close to the average of the initial conditions, thereby exhibiting behaviour similar to the simulation results, and, hence, validating our approach. It can further be observed that the error margin chosen is quite aggressive and if it was relaxed to 1% of the average value, the algorithm will converge substantially early.

## 6. CONCLUSION

In this article, the problem of consensus and average consensus is discussed in the presence of fixed delays in the network. Under the assumption of a strongly connected graph formed by the agents in the network, consensus is shown to be achieved asymptotically. A novel Max-Min Consensus based finite-time stopping criterion is introduced to distributively terminate the computation of consensus by the agents when each of them has reached within a pre-specified error bound. Further this algorithm is integrated with ratio consensus algorithm to prove the finite-time convergence to the average of initial conditions. The deviation achieved from the average of initial conditions by using the proposed distributed finite-time stopping criterion is also quantified. Furthermore, proper choice of the tolerance bound for stopping consensus algorithm can facilitate the convergence of average consensus arbitrarily close to the actual average of initial conditions. Furthermore, the practicality and real-time implementation of the proposed algorithm has been verified by testing it on a network of agents (Raspberry Pi devices) communicating over actual communication channels. The simulation and experimental results are in close agreement, thus, establishing the proposed method as a valid and practically attractive algorithm. This article is one of the very few attempts made in the direction of finite-time stopping of consensus and average consensus algorithms. To the best of the knowledge of the authors, the proposed algorithm is the simplest in terms of algorithmic and computational complexity.

Table II: Experimental results for the 5 node network using delay-free termination algorithm of [20].

| Agent   | Initial Value | Converged Value | Error  |
|---------|---------------|-----------------|--------|
| Agent 1 | 1000          | 225.05          | 174.95 |
| Agent 2 | 0             | 225.05          | 174.95 |
| Agent 3 | 200           | 225.05          | 174.95 |
| Agent 4 | 100           | 225.05          | 174.95 |
| Agent 5 | 700           | 225.05          | 174.95 |

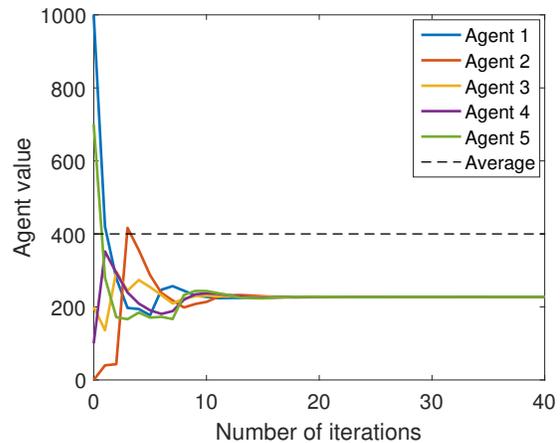

Figure 6: Experimental results for delay-free termination algorithm in presence of fixed delays on the network of 5 nodes shown in Figure 1(b).

Table III: Experimental results for 4 node ring network.

| Agent   | Initial Value | Converged Value | Error |
|---------|---------------|-----------------|-------|
| Agent 1 | 50            | 74.99           | 0.01  |
| Agent 2 | 70            | 74.99           | 0.01  |
| Agent 3 | 150           | 74.99           | 0.01  |
| Agent 4 | 30            | 74.99           | 0.01  |

Table IV: Experimental results for representative 5 node network.

| Agent   | Initial Value | Converged Value | Error |
|---------|---------------|-----------------|-------|
| Agent 1 | 1000          | 399.99          | 0.01  |
| Agent 2 | 0             | 399.99          | 0.01  |
| Agent 3 | 200           | 399.99          | 0.01  |
| Agent 4 | 100           | 399.99          | 0.01  |
| Agent 5 | 700           | 399.99          | 0.01  |

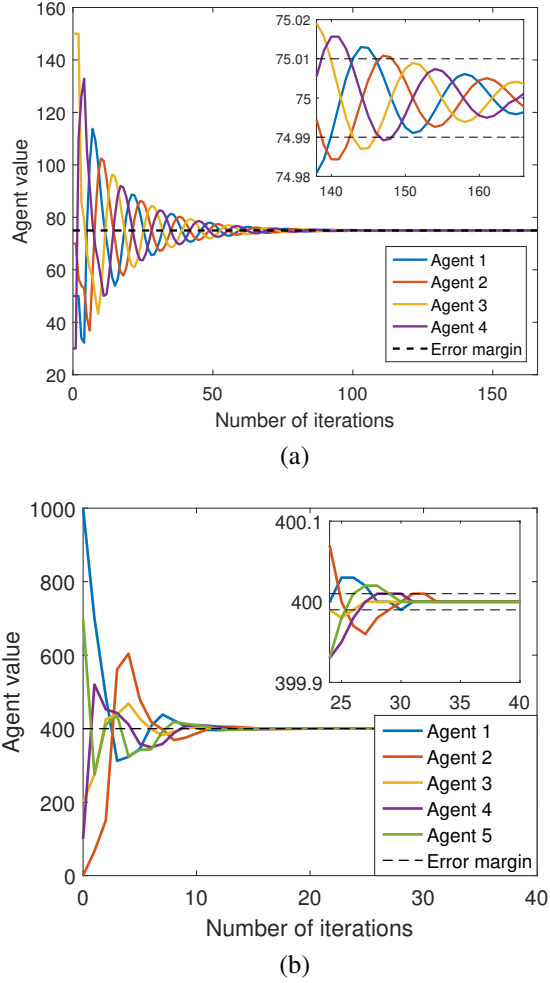

Figure 7: Experimental results for max-min based distributed finite-time criterion for average consensus protocol in presence of uniformly bounded fixed delays on the networks shown in Figure 1: (a) Ring network of 4 nodes (b) general representative network of 5 nodes.

# APPENDIX A
## PROOFS OF LEMMAS IN SECTION 3

### A. Proof of Lemma 3.1

*Proof.* From (1), we have,

$$x_i(k+1) = p_{ii}x_i(k) + \sum_{j \epsilon N_i^-} p_{ij}x_j(k-\tau_{ij}).$$

It follows that,

$$x_i(k+1) \leq p_{ii}x_{q_i}(k-\tau_{q_i}) + \sum_{j \epsilon N_i^-} p_{ij}x_{q_i}(k-\tau_{q_i})$$
$$= x_{q_i}(k-\tau_{q_i}), \text{ and}, \quad (19)$$
$$x_i(k+1) \geq p_{ii}x_{s_i}(k-\tau_{s_i}) + \sum_{j \epsilon N_i^-} p_{ij}x_{s_i}(k-\tau_{s_i})$$
$$= x_{s_i}(k-\tau_{s_i}). \quad (20)$$

By taking maximum over all nodes $i \in V$ in (19), it follows that,

$$x_i(k+1) \leq \max_{i \in V} x_i(k+1) \leq \max_{i \in V} x_{q_i}(k-\tau_{q_i}) = M(k), \quad (21)$$

for all $i \in V$. Similarly, by taking minimum over all nodes in (20), it follows that,

$$x_i(k+1) \geq \min_{i \in V} x_i(k+1) \geq \min_{i \in V} x_{s_i}(k-\tau_{s_i}) = m(k), \quad (22)$$

for all $i \in V$. When $k^{'} = k$ the proof for (6) and (7) follows from the definitions. For all time instants $k^{'} > k$, the proof for (6) and (7) is reached using strong induction and is presented below.

Define $\{q_m, \tau_{q_m}\} := \arg\max_{j \in V} x_{q_j}(k-\tau_{q_j})$ for some $q_m$ in $V, \tau_{q_m} \in \{0, 1, 2, ..., \bar{\tau}\}$. Thus, $x_{q_m}(k-\tau_{q_m})$ is the maximum value among all nodal states in the horizon of $\bar{\tau}$ into the past starting from iteration $k$, that is, by definition $x_{q_m}(k-\tau_{q_m}) = M(k)$.

Using (21), it follows that,

$$x_i(k+1) \leq \max_{j \in V} x_{q_j}(k-\tau_{q_j}) = x_{q_m}(k-\tau_{q_m}). \quad (23)$$

Suppose it is asserted that,

$$x_i(k^{'}) \leq \max_{j \in V} x_{q_j}(k-\tau_{q_j}) = x_{q_m}(k-\tau_{q_m}), \quad (24)$$
for $k' = \{k+2, ..., k+l\}$.

Thus, $x_i(k+2) \leq x_{q_m}(k-\tau_{q_m}), \cdots, x_i(k+l) \leq x_{q_m}(k-\tau_{q_m})$, for all $i = 1, 2, ..., n$. Define $\{q_m^{'}, \tau_{q_m^{'}}\} := \arg\max_{j \in V} x_{q_j}(k+l-\tau_{q_j})$ for some $q_m^{'}$ in $V, \tau_{q_m^{'}} \in \{0, 1, 2, ..., \bar{\tau}\}$, that is, $x_{q_m^{'}}(k-\tau_{q_m^{'}})$ is the maximum value among all nodal states in the horizon of $\bar{\tau}$ into the past starting from iteration $k+l$.

Using (21), it follows that $x_i(k+l+1) \leq \max_{j \in V} x_{q_j}(k+l-\tau_{q_j}) = x_{q_{m'}}(k+l-\tau_{q_{m'}})$, for all $i = 1, 2, ..., n$.

Consider the case when $l > \bar{\tau}$. Then $k+l-\tau_{q_{m'}} > k+\bar{\tau}-\bar{\tau} = k$. Thus, $k+l \geq k+l-\tau_{q_{m'}} \geq k+1$. Using (23) and (24), it follows that $x_i(k+l+1) \leq x_{q_{m'}}(k+l-\tau_{q_{m'}}) \leq x_{q_m}(k-\tau_{q_m})$.

Consider the case when $1 \leq l \leq \bar{\tau}$. Then, $k + l - \tau_{q_{m'}} \geq k + 1 - \bar{\tau}$. Furthermore, suppose $k \geq k + l - \tau_{q_{m'}}$. Then, $k \geq k + l - \tau_{q_{m'}} \geq k + 1 - \bar{\tau} > k - \bar{\tau}$. Thus, $x_{q_m}(k - \tau_{q_m}) = \max_{j \in V} x_{q_j}(k - \tau_{q_j}) = \max_{j \in V} \max_{\tau_{ij} \in \{0,1,2,...,\bar{\tau}\}} x_j(k - \tau_{ij}) \geq x_{q_{m'}}(k + l - \tau_{q_{m'}})$. It follows that,

$$x_i(k+l+1) \leq x_{q_{m'}}(k+l-\tau_{q_{m'}}) \leq x_{q_m}(k-\tau_{q_m}).$$

Now suppose $k < k + l - \tau_{q_{m'}}$. It follows that $k < k + l - \tau_{q_{m'}} \leq k+l$. Thus, from (21), $x_i(k+l+1) \leq x_{q_{m'}}(k+l-\tau_{q_{m'}})$. However, as the index $(k + l - \tau_{q_{m'}}) \in \{k+1, ..., k+l\}$, it follows from (23) and (24) that,

$$x_{q_{m'}}(k+l-\tau_{q_{m'}}) \leq x_{q_m}(k-\tau_{q_m}) = M(k).$$

It is thus established that, $x_i(k+l+1) \leq x_{q_{m'}}(k+l-\tau_{q_{m'}}) \leq x_{q_m}(k-\tau_{q_m}) = M(k)$.

This proves (6). Similarly, (7) can be proven. This completes the proof. □

### B. Proof of Lemma 3.2

*Proof.* By assumption for time instant $k' = k$, $x_i(k) < M(k)$. For all time instants $k' > k$, the proof is reached using strong induction. Note that $x_i(k+1) = p_{ii}x_i(k) + \sum_{j \in N_i^-} p_{ij}x_j(k-\tau_{ij})$. It follows that,

$$x_i(k+1) \leq p_{ii}x_i(k) + \sum_{j \in N_i^-} p_{ij} \max_{j \in V} x_{q_j}(k-\tau_{q_j}),$$

or, $x_i(k+1) \leq p_{ii}x_i(k) + \sum_{j \in N_i^-} p_{ij}M(k).$

Since, $x_i(k) < M(k)$, it implies that,

$$x_i(k+1) < p_{ii}M(k) + \sum_{j \in N_i^-} p_{ij}M(k),$$

or, $x_i(k+1) < M(k)$.

Suppose it is asserted that,

$$x_i(k+l') < M(k), \text{ for } l' = 2, ..., l. \quad (25)$$

It follows that,

$$x_i(k+l+1) = p_{ii}x_i(k+l) + \sum_{j \in N_i^-} p_{ij}x_j(k+l-\tau_{ij}).$$

Consider the case when $k \geq k+l-r$, where $r \in \{0,1,...,\bar{\tau}\}$. Then $k+2-\bar{\tau} \leq k+l-r \leq k$. As the index $(k+l-r) \in \{k+2-\bar{\tau},...,k\}$, it follows from the definition of $M(k)$ that,

$$x_j(k+l-r) \leq M(k).$$

Now consider the case when $k < k+l-r$. Using (25) and (6) it follows that,

$$x_j(k+l-r) \leq \max_{j \in V} x_{q_j}(k-\tau_{q_j}) = M(k).$$

Thus, $x_i(k+l+1) \leq p_{ii}x_i(k+l) + \sum_{j \in N_i^-} p_{ij}M(k)$. Hence,

$$x_i(k+l+1) < p_{ii}M(k) + \sum_{j \in N_i^-} p_{ij}M(k) = M(k),$$

implying, $x_i(k') < M(k)$, for all $k' \geq k+1$.

The proof for other inequality is similar to the proof above and is left to the reader. □

### C. Proof of Lemma 3.3

*Proof.* Since $\min x(k) < \max x(k)$, there exists a node $i$ such that $x_i(k) < \max x(k)$. Since $\max x(k) \leq \max_{j \in V} x_{q_j}(k - \tau_{q_j}) := M(k)$, it follows that

$$x_i(k) < M(k).$$

Consider any node $j$. By strong connectivity of $G$, there exists a directed path connecting nodes $i$ and $j$. Assume the shortest directed path connecting $i$ and $j$ is given by $(m_1, i)(m_2, m_1,) \cdots (j, m_{d_j-1})$. Using (1), it follows that

$$x_{m_1}(k+\bar{\tau}+1) = p_{m_1 i}x_i(k+\bar{\tau}-\tau_{m_1 i}) + \sum_{j \in N_{m_1}^- \cup \{m_1\} \setminus \{i\}} p_{m_1 j}x_j(k+\bar{\tau}-\tau_{m_1 j}).$$

As $(k + \bar{\tau} - \tau_{m_1 j}) \geq k$, using the definition of $M(k)$ and (6), it follows that

$$x_{m_1}(k+\bar{\tau}+1) \leq p_{m_1 i}x_i(k+\bar{\tau}-\tau_{m_1 i}) + \Sigma_{j \in N_{m_1}^- \cup \{m_1\} \setminus \{i\}} p_{m_1 j}M(k+\bar{\tau}).$$

However, since $M(k+\bar{\tau}) \leq M(k)$, it follows that,

$$x_{m_1}(k+\bar{\tau}+1) \leq p_{m_1 i}x_i(k+\bar{\tau}-\tau_{m_1 i}) + \Sigma_{j \in N_{m_1}^- \cup \{m_1\} \setminus \{i\}} p_{m_1 j}M(k).$$

As $k + \bar{\tau} - \tau_{m_1 i} \geq k$, using the fact that, $x_i(k) < M(k)$ and Lemma 3.2, it follows that

$$x_{m_1}(k+\bar{\tau}+1) < p_{m_1 i}M(k) + \sum_{j \in N_{m_1}^- \cup \{m_1\} \setminus \{i\}} p_{m_1 j}M(k).$$

Thus, $x_{m_1}(k+\bar{\tau}+1) < M(k)$. Thus, using Lemma 3.2, it follows that for all $k'' \geq k+\bar{\tau}+1$,

$$x_{m_1}(k'') < M(k).$$

Proceeding,

$$x_{m_2}(k+2\bar{\tau}+2) = p_{m_2 m_1}x_{m_1}(k+2\bar{\tau}+1-\tau_{m_2 m_1}) + \sum_{j \in N_{m_2}^- \cup \{m_2\} \setminus \{m_1\}} p_{m_2 j}x_j(k+2\bar{\tau}+1-\tau_{m_2 j}).$$

As $(k + 2\bar{\tau} + 1 - \tau_{m_2 j}) \geq (k + \bar{\tau} + 1) \geq k + 1$, using (6) it follows that

$$x_{m_2}(k+2\bar{\tau}+2) \leq p_{m_2m_1}x_{m_1}(k+2\bar{\tau}+1-\tau_{m_1m_2})+ \sum_{j\in N^-_{m_2}\cup\{m_2\}\setminus\{m_1\}} p_{m_2j}M(k).$$

Note that $(k+2\bar{\tau}+1-\tau_{m_1m_2}) \geq (k+\bar{\tau}+1)$ and $x_{m_1}(k+\bar{\tau}+1) < M(k)$. Using Lemma 3.2, it follows that

$$x_{m_2}(k+2\bar{\tau}+2) < p_{m_2m_1}M(k)+ \sum_{j\in N^-_{m_2}\cup\{m_2\}\setminus\{m_1\}} p_{m_2j}M(k).$$

Thus,

$$x_{m_2}(k+2\bar{\tau}+2) < M(k).$$

Thus, it can be concluded using Lemma 3.2, that for all $k'' \geq k+2\bar{\tau}+2$, $x_{m_2}(k'') < M(k)$. It follows that for all $k'' \geq k+d_j\bar{\tau}+d_j$, $x_j(k'') < M(k)$.

Note that, $k+D\bar{\tau}+D \geq k+d_j\bar{\tau}+d_j$, and, thus, for any $v \in V, v \neq i$, $k'' \geq k+D\bar{\tau}+D$ implies, $x_v(k'') < M(k)$.

It follows that, for all $k' \geq k+D(\bar{\tau}+1)$,

$$\max x(k') < \max_{j\in V} x_{q_j}(k-\tau_{q_j}) = M(k).$$

The proof for other inequality is similar to the proof above and is left to the reader. □

## APPENDIX B
## TIME-VARYING DELAY FRAMEWORK

The consensus model which incorporates time-varying delays is dealt in [12] and [13]. Under this model, node $i \in V$ updates its state at instant $k+1$ by:

$$x_i(k+1) = p_{ii}x_i(k) + \sum_{j\in N^-_i} p_{ij}x_j(k-\tau_{ij}(k)), \quad (26)$$

where $\tau_{ij}(k) \in \{0,1,2,...,\bar{\tau}\}$.

The assumptions A1, A2, A4 and A5 presented in Section 2 are valid for this case as well whereas assumption A3 is not valid as the link delays are time-varying but uniformly bounded by some finite integer $\bar{\tau}$. It is proven in [12] and [13] that the consensus update algorithm given by (26) converges asymptotically. To prove the applicability of finite-time termination algorithm proposed in Section 4 (*Algorithm 1*) for the time-varying delay case, the results presented in Section 3 are proved here considering update equation given by (26). Like the fixed delay case, let for node $i$, the maximum over all values held by its neighbors including itself currently and in the $\bar{\tau}$ past instants be denoted by $x_{q_i}(k-\tau_{q_i})$ and the minimum over all values held by its neighbors including itself currently and in the $\bar{\tau}$ past instants be denoted by $x_{s_i}(k-\tau_{s_i})$. The definitions of $q_i, \tau_{q_i}, s_i$ and $\tau_{s_i}$ are same as in (2) and (3).

Since, $\tau_{ij}(k) \in \{\tau_{ij}\}$, it follows that for $j \in N^-_i \cup \{i\}, x_j(k-\tau_{ij}(k)) \leq x_{q_i}(k-\tau_{q_i})$ and $x_j(k-r_{ij}(k)) \geq x_{s_i}(k-\tau_{s_i})$.

Furthermore, denote the maximum and minimum over all nodal values over the horizon in the past $\{k-\bar{\tau},...,k-1,k\}$ by $M(k)$ and $m(k)$ respectively as in (4) and (5).

**Lemma B.1.** *Consider the update rule (26). Then for all time instants $k' \geq k$ and for all $i \in V$,*

$$x_i(k') \leq \max_{j\in V} x_{q_j}(k-\tau_{q_j}) = M(k), \text{ and,} \quad (27)$$

$$x_i(k') \geq \min_{s\in V} x_{q_s}(k-\tau_{q_s}) = m(k). \quad (28)$$

*Proof.* From (26), we have

$$x_i(k+1) = p_{ii}x_i(k) + \sum_{j\epsilon N^-_i} p_{ij}x_j(k-\tau_{ij}(k)).$$

It follows that,

$$x_i(k+1) \leq p_{ii}x_{q_i}(k-\tau_{q_i}) + \sum_{j\epsilon N^-_i} p_{ij}x_{q_i}(k-\tau_{q_i})$$
$$= x_{q_i}(k-\tau_{q_i}), \text{ and,} \quad (29)$$
$$x_i(k+1) \geq p_{ii}x_{s_i}(k-\tau_{s_i}) + \sum_{j\epsilon N^-_i} p_{ij}x_{s_i}(k-\tau_{s_i})$$
$$= x_{s_i}(k-\tau_{s_i}). \quad (30)$$

By taking maximum over all nodes $i \in V$ in (29), it follows that,

$$x_i(k+1) \leq \max_{i\in V} x_i(k+1) \leq \max_{i\in V} x_{q_i}(k-\tau_{q_i}) = M(k), \quad (31)$$

for all $i \in V$. Similarly, by taking minimum over all nodes in (30), it follows that,

$$x_i(k+1) \geq \min_{i\in V} x_i(k+1) \geq \min_{i\epsilon V} x_{s_i}(k-\tau_{s_i}) = m(k), \quad (32)$$

for all $i \in V$. It should be noted that (31) and (32) are same as (21) and (22) respectively. Hence, the rest of the proof follows exactly the same way as the proof of *Lemma 3.1*.
□

**Lemma B.2.** *Consider a strongly connected graph $G = \{V,E\}$ running consensus protocol given by (26) with an initial condition $x(k)$. Let $i$ be a node such that $x_i(k) < M(k)$ and let $j$ be a node such that $x_j(k) > m(k)$, then for all time instants $k' \geq k$, $x_i(k') < M(k)$ and $x_j(k') > m(k)$.*

*Proof.* By assumption for time instant $k' = k, x_i(k) < M$. For all time instants $k' > k$, the proof is reached using strong induction. Note that

$$x_i(k+1) = p_{ii}x_i(k) + \sum_{j\in N^-_i} p_{ij}x_j(k-\tau_{ij}(k)).$$

It follows that

$$x_i(k+1) \leq p_{ii}x_i(k) + \sum_{j\in N^-_i} p_{ij}\max_{j\in V} x_{q_j}(k-\tau_{q_j}).$$

Then,

$$x_i(k+1) \leq p_{ii}x_i(k) + \sum_{j\in N^-_i} p_{ij}M(k).$$

Since, $x_i(k) < M(k)$, it follows that,

$$x_i(k+1) < p_{ii}M(k) + \sum_{j \in N_i^-} p_{ij}M(k).$$

Thus,

$$x_i(k+1) < M(k).$$

Suppose it is asserted that,

$$x_i(k+l') < M(k), \text{ for } l' = 2,...,l. \quad (33)$$

It follows that

$$x_i(k+l+1) = p_{ii}x_i(k+l) + \sum_{j \in N_i^-} p_{ij}x_j(k+l-\tau_{ij}(k+l)).$$

Consider the case when $k \geq (k+l-\tau_{ij}(k+l))$, then $k+2-\bar{\tau} \leq (k+l-\tau_{ij}(k+l)) \leq k$. As the index $(k+l-\tau_{ij}(k+l)) \in \{k+2-\bar{\tau},...,k\}$, it follows from the definition of $M(k)$ that

$$x_j(k+l-\tau_{ij}(k+l)) \leq M(k).$$

Now consider the case when $k < (k+l-\tau_{ij}(k+l))$. Using (33) and (27) it follows that

$$x_j(k+l-\tau_{ij}(k+l)) \leq \max_{j \in V} x_{q_j}(k-\tau_{q_j}) = M(k).$$

Thus, $x_i(k+l+1) \leq p_{ii}x_i(k+l) + \sum_{j \in N_i^-} p_{ij}M(k)$.
Hence, $x_i(k+l+1) < p_{ii}M(k) + \sum_{j \in N_i^-} p_{ij}M(k) = M(k)$.
This implies that,

$$x_i(k') < \max_{j \in V} x_{q_j}(k-\tau_{q_j}) = M(k) \text{ for all } k' \geq k+1.$$

The proof for other inequality is similar to the proof above and is left to the reader. □

Like in Section 3, consider the maximum and minimum value in the network, which is defined as, $\max x(k) := \max_{i \in V} x_i(k)$ and $\min x(k) := \min_{i \in V} x_i(k)$ respectively.

**Lemma B.3.** *Consider a strongly connected graph $G = \{V,E\}$ with an update rule for the consensus protocol given by (26) with an initial condition $x(k)$ such that $\min x(k) < \max x(k)$. Then for all $k' \geq k+D(\bar{\tau}+1)$, $\max x(k') < M(k)$ and $\min x(k') > m(k)$.*

*Proof.* Since $\min x(k) < \max x(k)$, there exists a node $i$ such that $x_i(k) < \max x(k)$. Since $\max x(k) \leq \max_{j \in V} x_{q_j}(k-\tau_{q_j}) := M(k)$, it follows that

$$x_i(k) < M(k).$$

Consider any node $j$. By strong connectivity of $G$, there exists a directed path connecting nodes $i$ and $j$. Assume the shortest directed path connecting $i$ and $j$ is given by $(m_1, i)(m_2, m_1)\cdots(j, m_{d_j-1})$. Using (26), it follows that

$$x_{m_1}(k+\bar{\tau}+1) = p_{m_1 i}x_i(k+\bar{\tau}-\tau_{m_1 i}(k+\bar{\tau}))+ \sum_{j \in N_{m_1}^- \cup \{m_1\}\setminus\{i\}} p_{m_1 j}x_j(k+\bar{\tau}-\tau_{m_1 j}(k+\bar{\tau})).$$

As $(k+\bar{\tau}-\tau_{m_1 j}(k+\bar{\tau})) \geq k$, using the definition of $M(k)$ and (27), it follows that

$$x_{m_1}(k+\bar{\tau}+1) \leq p_{m_1 i}x_i(k+\bar{\tau}-\tau_{m_1 j}(k+\bar{\tau}))+ \sum_{j \in N_{m_1}^- \cup \{m_1\}\setminus\{i\}} p_{m_1 j}M(k).$$

As $(k+\bar{\tau}-\tau_{m_1 j}(k+\bar{\tau})) \geq k$, using the fact that, $x_i(k) < M(k)$ and using Lemma B.2, it follows that

$$x_{m_1}(k+\bar{\tau}+1) < p_{m_1 i}M(k) + \sum_{j \in N_{m_1}^- \cup \{m_1\}\setminus\{i\}} p_{m_1 j}M(k).$$

Thus, $x_{m_1}(k+\bar{\tau}+1) < M(k)$. Thus, using Lemma B.2, it follows that for all $k'' \geq k+\bar{\tau}+1$,

$$x_{m_1}(k'') < M(k).$$

Proceeding further on the directed path,

$$x_{m_2}(k+2\bar{\tau}+2) = p_{m_2 m_1}x_{m_1}(k+2\bar{\tau}+1-\tau_{m_2 m_1}(k+2\bar{\tau}+1))+ \sum_{j \in N_{m_2}^- \cup \{m_2\}\setminus\{m_1\}} p_{m_2 j}x_j(k+2\bar{\tau}+1-\tau_{m_2 j}(k+2\bar{\tau}+1)).$$

As $(k+2\bar{\tau}+1-\tau_{m_2 j}(k+2\bar{\tau}+1)) \geq (k+\bar{\tau}+1) \geq k+1$, using (27), it follows that,

$$x_{m_2}(k+2\bar{\tau}+2) \leq p_{m_2 m_1}x_{m_1}(k+2\bar{\tau}+1-\tau_{m_2 m_1}(k+2\bar{\tau}+1)) + \sum_{j \in N_{m_2}^- \cup \{m_2\}\setminus\{m_1\}} p_{m_2 j}M(k).$$

Note that $(k+2\bar{\tau}+1-\tau_{m_2 m_1}(k+2\bar{\tau}+1)) \geq (k+\bar{\tau}+1)$ and $x_{m_1}(k+\bar{\tau}+1) < M(k)$. Using Lemma B.2, it follows that,

$$x_{m_2}(k+2\bar{\tau}+2) < p_{m_2 m_1}M(k) + \sum_{j \in N_{m_2}^- \cup \{m_2\}\setminus\{m_1\}} p_{m_2 j}M(k).$$

Thus,

$$x_{m_2}(k+2\bar{\tau}+2) < M(k).$$

Thus, it can be concluded using Lemma B.2, that for all $k'' \geq k+2\bar{\tau}+2$, $x_{m_2}(k'') < M(k)$. It follows that for all $k'' \geq k+d_j\bar{\tau}+d_j$, $x_j(k'') < M(k)$.
Note that, $k+D\bar{\tau}+D \geq k+d_j\bar{\tau}+d_j$, and, thus, for any $v \in V, v \neq i$, $k'' \geq k+D\bar{\tau}+D$ implies $x_v(k'') < M(k)$. It follows that, for all $k' \geq k+D(\bar{\tau}+1)$,

$$\max x(k') < \max_{j \in V} x_{q_j}(k-\tau_{q_j}) = M(k).$$

The proof for other inequality is similar to the proof above and is left to the reader. □

**Remark 3.** *It should be noted that the results of Lemma B.1, Lemma B.2 and Lemma B.3 are exactly the same as that of Lemma 3.1, Lemma 3.2 and Lemma 3.3 respectively and hence the results of Lemma 3.1, Lemma 3.2 and Lemma* 3.3 *hold for consensus with time-varying delays as well. Furthermore, the results of* Theorem 3.1*,* Corollary 3.1 *and* Corollary 3.2 *follow directly from the results of* Lemma 3.3*. Hence,* Theorem 3.1*,* Corollary 3.1 *and* Corollary 3.2 *hold true for consensus with time-varying delay framework.*

### A. Max and Min Consensus in Presence of Time-Varying Delays

Like the fixed delay case, the Maximum Consensus Protocol computes the maximum of the given initial node conditions $z(0) = [z_1(0)\ z_2(0)....z_n(0)]^T$ in a distributed manner. It takes $z(0)$ as an input and generates a sequence of node values based on the update rule for node $m$ given by (9) and (10). Similarly, the Minimum Consensus Protocol computes the minimum of the given initial node conditions $y(0) = [y_1(0)\ y_2(0)....y_n(0)]^T$ in a distributed manner. It takes $y(0)$ as an input and generates a sequence of node values based on the update rule for node $m$ given by (11) and (12).

In time-varying delay framework, it is possible that any in-neighbor $j$ of node $m$ can be sending multiple packets of $z_j$ to node $m$ while it is waiting for $\bar{\tau}$ iterations, but each packet received from $j$ during the waiting time contains the same information because node $j$ also updates its maximum consensus state only once in $\bar{\tau}$ iterations in accordance with (9) and (10). Similar logic holds for Min consensus protocol.

**Remark 4.** *It is worth noting that the results of* Lemma 3.4*,* Lemma 3.5*,* Lemma 3.6 *and* Lemma 3.7 *depend solely on Max-Min consensus update rules given by (9), (10), (11) and (12). Since the Max-Min consensus update rules are logically the same for time-varying delay case as well, the results of* Lemma 3.4*,* Lemma 3.5*,* Lemma 3.6 *and* Lemma 3.7 *hold true for time-varying delay framework. This entails that the results of associated* Lemma 3.8 *and* Theorem 3.2 *hold true for time-varying delay framework. This implies that the finite-time termination algorithm for consensus (*Algorithm 1*) proposed for the fixed delay case can be implemented without any modification even if the delays are time-varying and follow the update rule given by (26).*

## APPENDIX C
## QUANTIFYING DEVIATION FROM AVERAGE USING *Algorithm 2*

It has been shown in [24] that

$$\lim_{k \to \infty} \frac{x_i(k)}{w_i(k)} = c = \frac{\alpha}{\sigma} \quad (34)$$

Since, the two consensus protocols given by (13) and (14) will terminate when they reach within some specified bound $\rho > 0$, $x_i(k)$ and $w_i(k)$ do not converge to $\alpha$ and $\sigma$ respectively. Instead let $\tilde{\alpha}$ and $\tilde{\sigma}$ be the terminal values of $x_i$ and $w_i$ respectively upon termination of *Algorithm 2*, such that,

$$\alpha - \rho \leq \tilde{\alpha} \leq \alpha + \rho, \quad (35)$$

$$\sigma - \rho \leq \tilde{\sigma} \leq \sigma + \rho. \quad (36)$$

Using (34), (35) and (36), it follows that $\frac{\tilde{\alpha}-\rho}{\tilde{\sigma}+\rho} \leq c \leq \frac{\tilde{\alpha}+\rho}{\tilde{\sigma}-\rho}$. It is easy to see that $\frac{x_i}{w_i} - \frac{\tilde{\alpha}}{\tilde{\sigma}}\left(\frac{1-\frac{\rho}{\tilde{\alpha}}}{1+\frac{\rho}{\tilde{\sigma}}}\right) \geq \frac{x_i}{w_i} - c \geq \frac{x_i}{w_i} - \frac{\tilde{\alpha}}{\tilde{\sigma}}\left(\frac{1+\frac{\rho}{\tilde{\alpha}}}{1-\frac{\rho}{\tilde{\sigma}}}\right)$.
Thus, $\frac{\tilde{\alpha}}{\tilde{\sigma}} - \frac{\tilde{\alpha}}{\tilde{\sigma}}\left(\frac{1-\frac{\rho}{\tilde{\alpha}}}{1+\frac{\rho}{\tilde{\sigma}}}\right) \geq \frac{x_i}{w_i} - c \geq \frac{\tilde{\alpha}}{\tilde{\sigma}} - \frac{\tilde{\alpha}}{\tilde{\sigma}}\left(\frac{1+\frac{\rho}{\tilde{\alpha}}}{1-\frac{\rho}{\tilde{\sigma}}}\right)$.
Also, it is easy to see that as $\rho \to 0$, $\frac{x_i}{w_i} \to c$.

### ACKNOWLEDGMENT
The authors would like to thank the ARPA-E for supporting this research via ARPA-E Award No. DE-AR000071 for the project 'A Robust Distributed Framework for Flexible Power Grids'.

### REFERENCES

[1] M. Newman, A.-L. Barabasi, and D. J. Watts, *The structure and dynamics of networks*. Princeton University Press, 2006.

[2] F. Xia, L. T. Yang, L. Wang, and A. Vinel, "Internet of things," *International Journal of Communication Systems*, vol. 25, no. 9, p. 1101, 2012.

[3] M. Mesbahi and M. Egerstedt, *Graph theoretic methods in multiagent networks*. Princeton University Press, 2010.

[4] J. N. Tsitsiklis, "Problems in decentralized decision making and computation." DTIC Document, Tech. Rep., 1984.

[5] J. E. Boillat, "Load balancing and poisson equation in a graph," *Concurrency: Practice and Experience*, vol. 2, no. 4, pp. 289–313, 1990.

[6] A. Jadbabaie, J. Lin *et al.*, "Coordination of groups of mobile autonomous agents using nearest neighbor rules," *Automatic Control, IEEE Transactions on*, vol. 48, no. 6, pp. 988–1001, 2003.

[7] L. Xiao, S. Boyd, and S. Lall, "A scheme for robust distributed sensor fusion based on average consensus," in *Information Processing in Sensor Networks, 2005. IPSN 2005. Fourth International Symposium on*. IEEE, 2005, pp. 63–70.

[8] M. Andreasson, D. V. Dimarogonas, H. Sandberg, and K. H. Johansson, "Distributed control of networked dynamical systems: Static feedback, integral action and consensus," *Automatic Control, IEEE Transactions on*, vol. 59, no. 7, pp. 1750–1764, 2014.

[9] V. Blondel, J. M. Hendrickx, A. Olshevsky, J. Tsitsiklis *et al.*, "Convergence in multiagent coordination, consensus, and flocking," in *IEEE Conference on Decision and Control*, vol. 44, no. 3. IEEE; 1998, 2005, p. 2996.

[10] Z. Li and Z. Duan, *Cooperative control of multi-agent systems: a consensus region approach*. CRC Press, 2014.

[11] P.-A. Bliman, A. Nedic, and A. Ozdaglar, "Rate of convergence for consensus with delays," in *Decision and Control, 2008. CDC 2008. 47th IEEE Conference on*. IEEE, 2008, pp. 4849–4854.

[12] M. Cao, A. S. Morse, and B. D. Anderson, "Reaching a consensus in a dynamically changing environment: A graphical approach," *SIAM Journal on Control and Optimization*, vol. 47, no. 2, pp. 575–600, 2008.

[13] M. Cao, A. S. Morse, and B. Anderson, "Reaching an agreement using delayed information," in *Proceedings of the 45th IEEE Conference on Decision and Control*. IEEE, 2006, pp. 3375–3380.

[14] K. I. Tsianos and M. G. Rabbat, "The impact of communication delays on distributed consensus algorithms," *arXiv preprint arXiv:1207.5839*, 2012.

[15] ——, "Distributed consensus and optimization under communication delays," in *Communication, Control, and Computing (Allerton), 2011 49th Annual Allerton Conference on*. IEEE, 2011, pp. 974–982.

[16] A. D. Dominguez-Garcia and C. N. Hadjicostis, "Distributed algorithms for control of demand response and distributed energy resources," in *Decision and Control and European Control Conference (CDC-ECC), 2011 50th IEEE Conference on*. IEEE, 2011, pp. 27–32.

[17] T. Charalambous, Y. Yuan, T. Yang, W. Pan, C. N. Hadjicostis, and M. Johansson, "Distributed finite-time average consensus in digraphs in the presence of time delays," *IEEE Transactions on Control of Network Systems*, vol. 2, no. 4, pp. 370–381, 2015.

[18] J. He, P. Cheng, L. Shi, J. Chen, and Y. Sun, "Time synchronization in wsns: A maximum-value-based consensus approach," *IEEE Transactions on Automatic Control*, vol. 59, no. 3, pp. 660–675, 2014.

[19] Y. Zhang and S. Li, "From simplicity to complexity based on consensus: A case study," *arXiv preprint arXiv:1610.09482*, 2016.

[20] V. Yadav and M. V. Salapaka, "Distributed protocol for determining when averaging consensus is reached," in *45th Annual Allerton Conf*, 2007, pp. 715–720.

[21] R. Diestel, *Graph Theory*. Berlin, Germany: Springer-Verlag, 2006.
[22] R. A. Horn and C. R. Johnson, *Matrix analysis*. Cambridge university press, 2012.
[23] A. Nedić and A. Ozdaglar, "Convergence rate for consensus with delays," *Journal of Global Optimization*, vol. 47, no. 3, pp. 437–456, 2010.
[24] C. N. Hadjicostis and T. Charalambous, "Average consensus in the presence of delays in directed graph topologies," *IEEE Transactions on Automatic Control*, vol. 59, no. 3, pp. 763–768, 2014.
[25] J. Ding and A. Zhou, *Nonnegative matrices, positive operators, and applications*. World Scientific Singapore, 2009.
[26] E. Upton and G. Halfacree, *Raspberry Pi user guide*. John Wiley & Sons, 2014.
[27] E. Perahia and R. Stacey, *Next Generation Wireless LANS: 802.11 n and 802.11 ac*. Cambridge university press, 2013.
[28] R. Fielding, J. Gettys, J. Mogul, H. Frystyk, L. Masinter, P. Leach, and T. Berners-Lee, "Hypertext transfer protocol–http/1.1," Tech. Rep., 1999.
[29] K. R. Fall and W. R. Stevens, *TCP/IP illustrated, volume 1: The protocols*. addison-Wesley, 2011.
[30] M. Cantelon, M. Harter, T. Holowaychuk, and N. Rajlich, *Node. js in Action*. Manning, 2014.

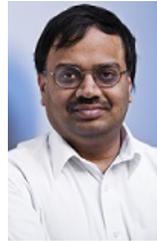

**Murti V. Salapka** received the B.Tech. degree in mechanical engineering from the Indian Institute of Technology, Madras, in 1991 and the M.S. and Ph.D. degrees in Mechanical Engineering from the University of California at Santa Barbara, in 1993 and 1997, respectively. He was a faculty member in the Electrical and Computer Engineering Department, Iowa State University, Ames, from 1997 to 2007. Currently, he is the Director of Graduate Studies and the Vincentine Hermes Luh Chair Professor in the Electrical and Computer Engineering Department, University of Minnesota, Minneapolis. His research interests include control and network science, nanoscience and single molecule physics. Dr. Salapka received the 1997 National Science Foundation CAREER Award.

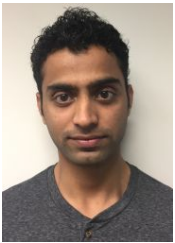

**Mangal Prakash** received the B.Tech degree in Electrical Engineering from National Institute of Technology, Durgapur, India and MS degree in Electrical Engineering from the University of Minnesota, USA. His research interests include control of network systems, probabilistic inference in graphical models and their applications to biological systems.

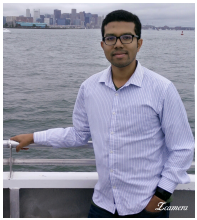

**Saurav Talukdar** received the B.Tech and M. Tech degree in Mechanical Engineering from Indian Institute of Technology, Bombay, India and is a PhD candidate in Mechanical Engineering at the University of Minnesota, Minneapolis, USA. His research interests include learning topology of dynamic systems, multi-agent systems and statistical physics.

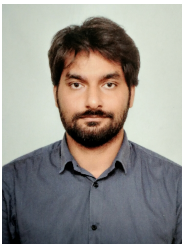

**Sandeep Attree** is a graduate student in electrical and computer engineering at the University of Minnesota. He obtained his B.Tech and M.Tech degrees from Indian Institute of Technology, Kanpur, India. His research interests include network optimization and distributed computing, with applications to smart and autonomous systems.

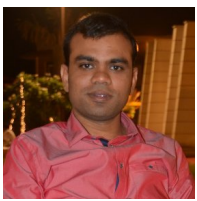

**Vikas Yadav** received the B.Tech. degree in electrical engineering from the Indian Institute of Technology, Kanpur, in 2000 and the M.S. and Ph.D. degrees in electrical engineering from Iowa State University, Ames, in 2007. He was with Garmin International Inc. from 2007 to 2013, Qualcomm Technologies Inc. from 2013 to 2015 and LG Electronics from 2015 to 2016. Presently he is a Principle Algorithm Architect at QuickLogic Corp., San Francisco, USA. His research interests include sensor fusion, distributed control design, self-organization and phase transition in large scale systems.